\documentclass{amsart}
\usepackage{amsopn}
\usepackage{amssymb, amscd}

\newcommand{\nc}{\newcommand}

\nc{\vg}{\mathfrak{v} } \nc{\wg}{\mathfrak{w} }
\nc{\zg}{\mathfrak{z} } \nc{\ngo}{\mathfrak{n} }
\nc{\kg}{\mathfrak{k} } \nc{\mg}{\mathfrak{m} }
\nc{\bg}{\mathfrak{b} } \nc{\ggo}{\mathfrak{g} }
\nc{\ggob}{\overline{\mathfrak{g}} } \nc{\sog}{\mathfrak{so} }
\nc{\sug}{\mathfrak{su} } \nc{\spg}{\mathfrak{sp} }
\nc{\slg}{\mathfrak{sl} } \nc{\glg}{\mathfrak{gl} }
\nc{\cg}{\mathfrak{c} } \nc{\rg}{\mathfrak{r} }
\nc{\hg}{\mathfrak{h} } \nc{\tg}{\mathfrak{t} }
\nc{\ug}{\mathfrak{u} } \nc{\dg}{\mathfrak{d} }
\nc{\ag}{\mathfrak{a} } \nc{\pg}{\mathfrak{p} }
\nc{\sg}{\mathfrak{s} } \nc{\pca}{\mathcal{P}}
\nc{\nca}{\mathcal{N}} \nc{\lca}{\mathcal{L}} \nc{\oca}{\mathcal{O}}
\nc{\mca}{\mathcal{M}} \nc{\tca}{\mathcal{T}} \nc{\aca}{\mathcal{A}}
\nc{\cca}{\mathcal{C}} \nc{\sca}{\mathcal{S}}

\nc{\vp}{\varphi} \nc{\ddt}{{\small \frac{{\rm d}}{{\rm d}t}}}
\nc{\im}{\mathtt{i}}

\nc{\SO}{{\mathrm SO}} \nc{\Spe}{{\mathrm Sp}} \nc{\Sl}{{\mathrm
SL}} \nc{\SU}{{\mathrm SU}} \nc{\Or}{{\mathrm O}} \nc{\U}{{\mathrm
U}} \nc{\Gl}{{\mathrm GL}} \nc{\Se}{{\mathrm S}} \nc{\Cl}{{\mathrm
Cl}} \nc{\Spein}{{\mathrm Spin}} \nc{\Pin}{{\mathrm Pin}}

\nc{\RR}{{\Bbb R}} \nc{\HH}{{\Bbb H}} \nc{\CC}{{\Bbb C}}
\nc{\ZZ}{{\Bbb Z}} \nc{\FF}{{\Bbb F}} \nc{\NN}{{\Bbb N}}
\nc{\QQ}{{\Bbb Q}} \nc{\PP}{{\Bbb P}}

\nc{\vs}{\vspace{.2cm}} \nc{\vsp}{\vspace{1cm}}
\nc{\ip}{\langle\cdot,\cdot\rangle} \nc{\la}{\langle}
\nc{\ra}{\rangle} \nc{\unm}{\frac{1}{2}} \nc{\unc}{\frac{1}{4}}
\nc{\und}{\frac{1}{16}} \nc{\no}{\vs\noindent}
\nc{\lam}{\Lambda^2\ngo^*\otimes\ngo} \nc{\tangz}{{\rm T}^{\rm Zar}}
\nc{\nor}{{\sf n}} \nc{\eigen}{(k_1<...<k_r;d_1,...,d_r)}
\nc{\eigencero}{(0<k_2<...<k_r;d_1,...,d_r)} \nc{\mum}{/\!\!/}
\nc{\kir}{/\!\!/\!\!/}

\nc{\He}{\operatorname{Hess}} \nc{\ad}{\operatorname{ad}}
\nc{\Ad}{\operatorname{Ad}} \nc{\rank}{\operatorname{rank}}
\nc{\Irr}{\operatorname{Irr}} \nc{\End}{\operatorname{End}}
\nc{\Aut}{\operatorname{Aut}} \nc{\Inn}{\operatorname{Inn}}
\nc{\Der}{\operatorname{Der}} \nc{\Ker}{\operatorname{Ker}}
\nc{\Iso}{\operatorname{I}} \nc{\Diff}{\operatorname{D}}
\nc{\Lie}{\operatorname{L}} \nc{\tr}{\operatorname{tr}}
\nc{\dif}{\operatorname{d}} \nc{\sen}{\operatorname{sen}}
\nc{\modu}{\operatorname{mod}} \nc{\Ric}{\operatorname{Ric}}
\nc{\Ricac}{\operatorname{Ric^{ac}}}
\nc{\Ricg}{\operatorname{Ric^{\gamma}}}
\nc{\Ricc}{\operatorname{Ric^{c}}} \nc{\sym}{\operatorname{sym}}
\nc{\symac}{\operatorname{sym^{ac}}}
\nc{\symc}{\operatorname{sym^{c}}} \nc{\scalar}{\operatorname{sc}}
\nc{\grad}{\operatorname{grad}} \nc{\ricci}{\operatorname{ric}}
\nc{\ricciac}{\operatorname{ric^{ac}}}
\nc{\riccic}{\operatorname{ric^{c}}}
\nc{\riccig}{\operatorname{ric^{\gamma}}}
\nc{\Rin}{\operatorname{M}} \nc{\Le}{\operatorname{L}}
\nc{\tang}{\operatorname{T}} \nc{\level}{\operatorname{level}}
\nc{\rad}{\operatorname{r}} \nc{\abel}{\operatorname{ab}}

\newtheorem{theorem}{Theorem}[section]
\newtheorem{proposition}[theorem]{Proposition}
\newtheorem{corollary}[theorem]{Corollary}

\newtheorem{definition}[theorem]{Definition}
\newtheorem{remark}[theorem]{Remark}

\newtheorem{example}[theorem]{Example}

\title{Minimal metrics on nilmanifolds}

\author{Jorge Lauret}

\address{FaMAF and CIEM, Universidad Nacional de C\'ordoba, 5000 C\'ordoba, Argentina}
\email{lauret@mate.uncor.edu}

\thanks{
Supported by CONICET and Guggenheim Foundation fellowships, and a
grant from FONCyT (Argentina).}

\begin{document}

\maketitle

\begin{abstract}
A left invariant metric on a nilpotent Lie group is called {\it
minimal}, if it minimizes the norm of the Ricci tensor among all
left invariant metrics with the same scalar curvature.  Such metrics
are unique up to isometry and scaling and the groups admitting a
minimal metric are precisely the nilradicals of (standard) Einstein
solvmanifolds.  If $N$ is endowed with an invariant symplectic,
complex or hypercomplex structure, then minimal compatible metrics
are also unique up to isometry and scaling.  The aim of this paper
is to give more evidence of the existence of minimal metrics, by
presenting several explicit examples.  This also provides many
continuous families of symplectic, complex and hypercomplex
nilpotent Lie groups.  A list of all known examples of Einstein
solvmanifolds is also given.

\end{abstract}

\section{Introduction}\label{intro}

A nilpotent Lie group $N$ can never admit an Einstein left invariant
metric, unless it is abelian.  A way of getting as close as possible
to this would be by defining a left invariant metric $\ip$ on $N$ to
be {\it minimal} if
$$
||\ricci_{\ip}||=\min \left\{ ||\ricci_{\ip'}|| :
\scalar(\ip')=\scalar(\ip)\right\},
$$
where $\ip'$ runs over all left invariant metrics on $N$ and
$\ricci_{\ip}$, $\scalar(\ip)$ denote the Ricci tensor and the
scalar curvature , respectively.  Indeed,
$$
||\ricci_{\ip}-\scalar(\ip)/n\ip||^2=||\ricci_{\ip}||^2
-\scalar(\ip)^2/n.
$$
A left invariant metric is always identified with the corresponding
inner product $\ip$ on the Lie algebra $\ngo$ of $N$.  The following
conditions on $\ip$ are equivalent to minimality and show that such
metrics are special from many other points of view:
\begin{itemize}
\item[(i)] $\ip$ is a Ricci soliton metric: the solution $\ip_t$
with initial point $\ip_0=\ip$ to the normalized Ricci flow
$$
\ddt\ip_t=-2\ricci_{\ip_t}-2||\ricci_{\ip_t}||^2\ip_t,
$$
(under which $\scalar(\ip_t)$ is constant in time) remains isometric
to $\ip$, that is, $\ip_t=\vp_t^*\ip$ for some one parameter group
of diffeomorphisms $\{\vp_t\}$ of $N$ (see \cite{soliton}).

\item[(ii)] $(N,\ip)$ admits a standard metric solvable extension
$(S ,\ip')$ which is Einstein: the Lie algebra $\sg$ of $S$ is given
by the orthogonal decomposition $\sg=\ag\oplus\ngo$, where $\ag$ is
abelian, $\ngo=[\sg,\sg]$ and $\ip'|{\ngo\times\ngo}=\ip$ (see
Section \ref{einstein}).

\item[(iii)] $\ip$ is a quasi-Einstein metric:
$$
\ricci_{\ip}=c\ip+L_X\ip
$$
for some $C^{\infty}$ vector field $X$ in $N$ and $c\in\RR$, where
$L_X\ip$ denotes the usual Lie derivative (see \cite{soliton}).

\item[(iv)] $\Ric_{\ip}=cI+D$ for some $c\in\RR$ and
$D\in\Der(\ngo)$, where $\Ric_{\ip}$ is the Ricci operator and
$\Der(\ngo)$ is the space of all derivations of $\ngo$.
\end{itemize}

The uniqueness up to isometry and scaling of minimal metrics on a
given nilpotent Lie group was proved in \cite{soliton}.  The
existence question is still nebulous; the only known obstruction
until now is that $\ngo$ has to admit an $\NN$-gradation, that is, a
direct sum decomposition $\ngo=\ngo_{k_1}\oplus...\oplus\ngo_{k_r}$,
$k_i\in\NN$, such that
$[\ngo_{k_i},\ngo_{k_j}]\subset\ngo_{k_i+k_j}$.  Such gradation is
defined by the symmetric derivation $D$ from condition (iv) above,
which was proved to have natural eigenvalues (up to scaling) by J.
Heber \cite{Hbr} in the context of Einstein solvmanifolds.
\begin{itemize}
\item[] {\bf Problem.} Does every $\NN$-graded nilpotent Lie algebra
admit a minimal metric?.
\end{itemize}
In view of equivalence (ii) above, what we are wondering is if any
$\NN$-graded nilpotent Lie group can be the nilradical of a standard
Einstein solvmanifold.  It is perhaps too optimistic to expect an
affirmative answer to this question, even in the two-step nilpotent
case.  The interplay with Einstein solvmanifolds provides a lot of
examples of nilpotent Lie groups admitting minimal metrics.  The
rather long list of known examples given in Section \ref{einstein}
shows that although the answer to the above question might be no,
there is a great deal of nilpotent Lie groups admitting a minimal
metric and also, most of those which are distinguished in some way
do so.

\vs

The search for a canonical metric also makes sense when there is a
given invariant geometric structure on $N$, as for example a
symplectic, complex or hypercomplex structure.  With these
structures in mind, we can define an {invariant geometric structure}
as a tensor on $N$ defined by left translation of a tensor $\gamma$
on $\ngo$ (or a set of tensors), usually non-degenerate in some way,
which satisfies a suitable integrability condition
\begin{equation}\label{closedgint}
{\rm IC}(\gamma,\mu)=0,
\end{equation}
involving only $\gamma$ and the Lie bracket $\mu$ of $\ngo$.  This
will allow us to study these three classes of structures and maybe
some other ones of similar characteristics at the same time.

The pair $(N,\gamma)$ will be called a {\it class-$\gamma$ nilpotent
Lie group}, and $N$ will be assumed to be simply connected for
simplicity.  A left invariant Riemannian metric on $N$ is said to be
{\it compatible} with $(N,\gamma)$ if the corresponding inner
product $\ip$ on $\ngo$ satisfies an orthogonality condition
\begin{equation}\label{ortconint}
{\rm OC}(\gamma,\ip)=0,
\end{equation}
in which only $\ip$ and $\gamma$ are involved.  We denote by
$\cca=\cca(N,\gamma)$ the set of all left invariant metrics on $N$
which are compatible with $(N,\gamma)$.  The pair $(\gamma,\ip)$
with $\ip\in\cca$ will often be referred to as a {\it class-$\gamma$
metric structure}.

The Ricci tensor has always been a very useful tool to deal with the
existence of distinguished metrics, and since the geometric
structure under consideration should be involved in the definition
of such a metric, we consider the {\it invariant Ricci operator}
$\Ricg_{\ip}$ (and the {\it invariant Ricci tensor}
$\riccig_{\ip}=\la\Ricg_{\ip}\cdot,\cdot\ra$), that is, the
orthogonal projection of the Ricci operator $\Ric_{\ip}$ onto the
subspace of those symmetric maps of $\ngo$ leaving $\gamma$
invariant. D. Blair, S. Ianus and A. Ledger \cite{BlrIns, BlrLdg,
Blr} have proved in the compact case that metrics satisfying
\begin{equation}\label{unconditionint}
\riccig_{\ip}=0
\end{equation}
are very special in symplectic (so called metrics with hermitian
Ricci tensor) and contact geometry, as they are precisely the
critical points of two very natural curvature functionals on $\cca$:
the total scalar curvature functional $S$ and a functional $K$
measuring how far are the metrics of being K$\ddot{{\rm a}}$hler or
Sasakian, respectively (see Section \ref{symp} for further
information).

Unfortunately, for a non-abelian nilpotent Lie group, condition
(\ref{unconditionint}) is too strong for the classes of structures
we have in mind, and hence it is natural to try to get as close as
possible to this unattainable goal.  In this light, a metric
$\ip\in\cca(N,\gamma)$ is called {\it minimal} if it minimizes the
functional $||\riccig_{\ip}||^2=\tr(\Ricg_{\ip})^2$ on the set of
all compatible metrics with the same scalar curvature.  Recall that
for $\gamma=0$ (i.e. when we are not considering any structure) the
Ricci invariant tensor coincides with the usual Ricci tensor and so
we get precisely minimal metrics as defined at the beginning of this
section.

We may also try to improve the metric via the evolution flow
\begin{equation}\label{flow}
\ddt \ip_t=\pm\riccig_{\ip_t},
\end{equation}
whose fixed points are precisely metrics satisfying
(\ref{unconditionint}) (the choice of the right sign depend on the
class of structure). In the symplectic case, this flow is called the
anticomplexified Ricci flow and has been recently studied by H-V Le
and G. Wang \cite{LeWng}. Of particular significance are then those
metrics for which the solution to the normalized flow (under which
the scalar curvature is constant in time) remains isometric in time
to the initial metric.  Such special metrics will be called {\it
invariant Ricci solitons}.  The following theorem was obtained by
using deep results from geometric invariant theory concerning the
moment map for a linear action of a reductive Lie group (see Section
\ref{var}).

\begin{theorem}\label{equiv2gint}\cite{minimal}
Let $(N,\gamma)$ be a nilpotent Lie group endowed with an invariant
geometric structure $\gamma$ (non-necessarily integrable). Then the
following conditions on a left invariant Riemannian metric $\ip$
which is compatible with $(N,\gamma)$ are equivalent:
\begin{itemize}
\item[(i)] $\ip$ is minimal.

\item[(ii)] $\ip$ is an invariant Ricci soliton.

\item[(iii)] $\Ricg_{\ip}=cI+D$ for some $c\in\RR$, $D\in\Der(\ngo)$.
\end{itemize}
Moreover, there is at most one compatible left invariant metric on
$(N,\gamma)$ up to isometry  (and scaling) satisfying any of the
above conditions.
\end{theorem}

\begin{corollary}\label{isoiso}\cite{minimal}
Let $\gamma,\gamma'$ be two geometric structures on a nilpotent Lie
group $N$, and assume that they admit minimal compatible metrics
$\ip$ and $\ip'$, respectively.  Then $\gamma$ is isomorphic to
$\gamma'$ if and only if there exists $\vp\in\Aut(\ngo)$ and $c>0$
such that $\gamma'=\vp.\gamma$ and
$$
\la\vp X,\vp Y\ra'=c\la X,Y\ra \qquad \forall X,Y\in\ngo.
$$
In particular, if $\gamma$ and $\gamma'$ are isomorphic then their
respective minimal compatible metrics are necessarily isometric up
to scaling (recall that $c=1$ when $\scalar(\ip)=\scalar(\ip')$).
\end{corollary}

A major obstacle to classify geometric structures is the lack of
invariants.  The uniqueness result in the above theorem and its
corollary gives rise to a useful tool to distinguish two geometric
structures; indeed, if they are isomorphic then their respective
minimal compatible metrics (if any) have to be isometric.  One
therefore can eventually distinguish geometric structures with
Riemannian data, which suddenly provides us with a great deal of
invariants.  This will be used in this paper to give explicit
continuous families of pairwise non-isomorphic geometric structures
in low dimensions, mainly by using only one Riemannian invariant:
the eigenvalues of the Ricci operator.  To actually find the
candidates for such families we apply a variational method which is
explained in Section \ref{var}.

Existence of minimal compatible metrics is proved for all
$4$-dimensional symplectic structures and a curve in dimension $6$,
two curves of abelian complex structures on the Iwasawa manifold and
several continuous families depending on various parameters of
abelian and non-abelian hypercomplex structures in dimension $8$.
It is finally showed in Section \ref{einstein} that if one considers
no structure (i.e. $\gamma=0$), then the `moment map' approach
proposed in \cite{minimal} can be also applied to the study of
Einstein solvmanifolds, obtaining many of the uniqueness and
structure results proved by J. Heber in \cite{Hbr}.

The existence problem is also far to be solved in this case; the
theorem does not even suggest when such a distinguished metric does
exist. How special are the symplectic or complex structures
admitting a minimal metric?. So far, we know how to deal with this
`existence question' only by giving several explicit examples, which
is the aim of this paper. The neat `algebraic' characterization
(iii) in Theorem \ref{equiv2gint} will be very useful. It turns out
that in low dimensions the structures in general tend to admit a
minimal compatible metric. At the moment, the only counterexamples
we have to the existence question are the characteristically
nilpotent Lie algebras (i.e. $\Der(\ngo)$ is nilpotent) admitting a
symplectic structure recently found by D. Burde in \cite{Brd}.

\begin{remark}{\rm
In \cite{classi}, by taking advantage of the interplay with
invariant theory, we describe the moduli space of all isomorphism
classes of geometric structures on nilpotent Lie groups of a given
class and dimension admitting a minimal compatible metric, as the
disjoint union of semi-algebraic varieties which are homeomorphic to
categorical quotients of suitable linear actions of reductive Lie
groups. Such special geometric structures can therefore be
distinguished by using invariant polynomials.  }
\end{remark}

\section{Variety of compatible metrics}\label{var}

Let us consider as a parameter space for the set of all real
nilpotent Lie algebras of a given dimension $n$, the set $\nca$ of
all nilpotent Lie brackets on a fixed $n$-dimensional real vector
space $\ngo$.  If
$$
V=\lam=\{\mu:\ngo\times\ngo\mapsto\ngo : \mu\; \mbox{skew-symmetric
bilinear map}\},
$$
then
$$
\nca=\{\mu\in V:\mu\;\mbox{satisfies Jacobi and is nilpotent}\}
$$
is an algebraic subset of $V$.  Indeed, the Jacobi identity and the
nilpotency condition are both determined by zeroes of polynomials.

We fix a tensor $\gamma$ on $\ngo$ (or a set of tensors), and let
$G_{\gamma}$ denote the subgroup of $\Gl(n)$ preserving $\gamma$.
These groups act naturally on $V$ by
\begin{equation}\label{action}
g.\mu(X,Y)=g\mu(g^{-1}X,g^{-1}Y), \qquad X,Y\in\ngo, \quad
g\in\Gl(n),\quad \mu\in V,
\end{equation}
and leave $\nca$ invariant.  Consider the subset
$\nca_{\gamma}\subset\nca$ given by
$$
\nca_{\gamma}=\{\mu\in\nca:{\rm IC}(\gamma,\mu)=0\},
$$
that is, those nilpotent Lie brackets for which $\gamma$ is
integrable (see (\ref{closedgint})). $\nca_{\gamma}$ is also an
algebraic variety since ${\rm IC}(\gamma,\mu)$ is always linear on
$\mu$ (at least in the cases we have in mind: symplectic, complex
and hypercomplex). Recall that
$$
W_{\gamma}=\{\mu\in V:{\rm IC}(\gamma,\mu)=0\}
$$
is a $G_{\gamma}$-invariant linear subspace of $V$, and
$\nca_{\gamma}=\nca\cap W_{\gamma}$.

For each $\mu\in\nca$, let $N_{\mu}$ denote the simply connected
nilpotent Lie group with Lie algebra $(\ngo,\mu)$.  Fix an inner
product $\ip$ on $\ngo$ compatible with $\gamma$, that is, such that
(\ref{ortconint}) holds. We identify each $\mu\in\nca_{\gamma}$ with
a class-$\gamma$ metric structure
\begin{equation}\label{ideg}
\mu\longleftrightarrow   (N_{\mu},\gamma,\ip),
\end{equation}
where all the structures are defined by left invariant translation.
Therefore, each $\mu\in\nca_{\gamma}$ can be viewed in this way as a
metric compatible with the class-$\gamma$ nilpotent Lie group
$(N_{\mu},\gamma)$, and two metrics $\mu,\lambda$ are compatible
with the same geometric structure if and only if they live in the
same $G_{\gamma}$-orbit. Indeed, the action of $G_{\gamma}$ on
$\nca_{\gamma}$ has the following interpretation: each $\vp\in
G_{\gamma}$ determines a Riemannian isometry preserving the
geometric structure
$$
(N_{\vp.\mu},\gamma,\ip)\mapsto
(N_{\mu},\gamma,\la\vp\cdot,\vp\cdot\ra)
$$
by exponentiating the Lie algebra isomorphism
$\vp^{-1}:(\ngo,\vp.\mu)\mapsto(\ngo,\mu)$.  We then have the
identification $G_{\gamma}.\mu=\cca(N_{\mu},\gamma)$, and more in
general the following

\begin{proposition}\label{upto}\cite{minimal}
Every class-$\gamma$ metric structure $(N',\gamma',\ip')$ on a
nilpotent Lie group $N'$ of dimension $n$ is isometric-isomorphic to
a $\mu\in\nca_{\gamma}$.
\end{proposition}

According to the above proposition and identification (\ref{ideg}),
the orbit $G_{\gamma}.\mu$ parameterizes all the left invariant
metrics which are compatible with $(N_{\mu},\gamma)$ and hence we
may view $\nca_{\gamma}$ as the space of all class-$\gamma$ metric
structures on nilpotent Lie groups of dimension $n$. Since two
metrics $\mu,\lambda\in\nca_{\gamma}$ are isometric if and only if
they live in the same $K_{\gamma}$-orbit, where
$K_{\gamma}=G_{\gamma}\cap\Or(\ngo,\ip)$ (see
\cite[Appendix]{minimal}), we have that $\nca_{\gamma}/K_{\gamma}$
parameterizes class-$\gamma$ metric nilpotent Lie groups of
dimension $n$ up to isometry and $G_{\gamma}.\mu/K_{\gamma}$ do the
same for all the compatible metrics on $(N_{\mu},\gamma)$.

In the search for the best compatible metric, it is natural to
consider the functional $F:\nca_{\gamma}\mapsto\RR$ given by
$F(\mu)=\tr(\Ricg_{\mu})^2$, which in some sense measures how far
the metric $\mu$ is from satisfying (\ref{unconditionint}). The
critical points of $F/||\mu||^4$ on the projective algebraic variety
$\PP\nca_{\gamma}\subset\PP V$ (which is equivalent to normalize by
the scalar curvature since $\scalar(\mu)=-\unc ||\mu||^2$), may
therefore be considered compatible metrics of particular
significance.

A crucial fact of this approach is that the moment map
$m_{\gamma}:V\mapsto\pg_{\gamma}$ for the action of $G_{\gamma}$ on
$V$ is proved to be
$$
m_{\gamma}(\mu)=8\Ricg_{\mu}, \qquad \forall\;\mu\in\nca_{\gamma},
$$
where $\pg_{\gamma}$ is the space of symmetric maps of $(\ngo,\ip)$
leaving $\gamma$ invariant (i.e.
$\ggo_{\gamma}=\kg_{\gamma}\oplus\pg_{\gamma}$ is a Cartan
decomposition).  This allows us to use strong and well-known results
on the moment map due to F. Kirwan \cite{Krw1} and L. Ness
\cite{Nss}, and proved by A. Marian \cite{Mrn} in the real case (we
also refer to \cite[Section 3]{minimal} and \cite{strata} for
further information). Indeed, since $F$ becomes a scalar multiple of
the square norm of the moment map, we obtain the following

\begin{theorem}\label{equiv1gint}\cite{Mrn}
Let $F:\PP \nca_{\gamma}\mapsto\RR$ be defined by
$F([\mu])=\tr(\Ricg_{\mu})^2/||\mu||^4$.  Then for $\mu\in
\nca_{\gamma}$ the following conditions are equivalent:
\begin{itemize}
\item[(i)] $[\mu]$ is a critical point of $F$.

\item[(ii)] $F|_{G_{\gamma}.[\mu]}$ attains its minimum value at $[\mu]$.

\item[(iii)] $\Ricg_{\mu}=cI+D$ for some $c\in\RR$, $D\in\Der(\mu)$.
\end{itemize}
Moreover, all the other critical points of $F$ in the orbit
$G_{\gamma}.[\mu]$ lie in $K_{\gamma}.[\mu]$.
\end{theorem}

The equivalence between (i) and (iii) in Theorem \ref{equiv2gint},
as well as the uniqueness result, follow then almost directly from
the above theorem.  We note that Theorem \ref{equiv1gint} also gives
a variational method to find minimal compatible metrics, by
characterizing them as the critical points of a natural curvature
functional (see \cite[Example 5.2]{minimal} for an explicit
application).

Most of the results obtained in this paper are still valid for
general Lie groups, although some considerations have to be
carefully taken into account (see \cite[Remark 4.6]{minimal}).

\section{Symplectic structures}\label{symp}

Let $(M,\omega)$ be a symplectic manifold, that is, a differentiable
manifold $M$ endowed with a global $2$-form $\omega$ which is closed
($\dif\omega=0$) and non-degenerate ($\omega^n\ne 0$).  A Riemannian
metric $g$ on $M$ is said to be {\it compatible} with $\omega$ if
there exists an almost-complex structure $J_g$ (i.e. a
$(1,1)$-tensor field with $J_g^2=-I$) such that
$$
\omega=g(\cdot,J_g\cdot).
$$
In that case $J_g$ is uniquely determined by $g$, and one may also
define that an almost-complex structure $J$ is {\it compatible} with
$\omega$ if
$$
g_J=\omega(\cdot,J\cdot)
$$
determines a Riemannian metric, which is again uniquely determined
by $J$.  In such a way we are really talking about compatible pairs
$(g,J)=(g,J_g)=(g_J,J)$, and the triple $(\omega,g,J_g)$ is called
an {\it almost-K$\ddot{{\rm a}}$hler} structure on $M$.

It is well known that for any symplectic manifold there always exist
a compatible metric.  Moreover, the space $\cca=\cca(M,\omega)$ of
all compatible metrics is usually huge; recall for instance that the
group of all symplectomorphisms (i.e. diffeomorphisms $\vp$ of $M$
such that $\vp^*\omega= \omega$) acts on $\cca$.

We fix from now on a symplectic manifold $(M,\omega)$.  Let $\Ric_g$
and $\nabla_g$ denote the Ricci operator and the Levi-Civita
conexion of a compatible metric $g$, respectively.  The most famous
conditions to ask g to satisfy are Einstein (i.e. $\Ric_g=cI$) and
K$\ddot{{\rm a}}$hler (i.e. $\nabla_gJ_g =0$), which are both very
strong and share the following property.

\begin{definition}\label{hrt}
{\rm We say that $g$ has {\it hermitian} (or $J$-{\it invariant})
{\it Ricci tensor} or that $J_g$ is {\it harmonic}, if $\Ric_gJ_g=
J_g\Ric_g$.}
\end{definition}
Examples of compatible metrics with hermitian Ricci tensor which are
neither Einstein nor K$\ddot{{\rm a}}$hler are known in any
dimension $2n\geq 6$ (see \cite{ApsGdc, DvdMsk,LeWng}).  It is
proved in \cite{minimal} that a symplectic nilpotent Lie group can
never admit a compatible metric with hermitian Ricci tensor unless
it is abelian.

A classical approach to searching for distinguished metrics is the
variational one, that is, to consider critical points of natural
functionals of the curvature on the space of all metrics of a given
class.  For instance, if $M$ is compact, D. Hilbert \cite{Hlb}
proved that Einstein metrics on $M$ are precisely the critical
points of the total scalar curvature functional
$$
S:\mca_1\mapsto\RR, \qquad S(g)=\int_M\scalar(g)\dif\nu_g,
$$
where $\mca_1$ is the space of all Riemannian metrics on $M$ with
volume equal to $1$.  Since the set of compatible metrics $\cca$ is
smaller, one should expect a weaker critical point condition for
$S:\cca\mapsto\RR$.  Another natural functional in our setup would
be
$$
K:\cca\mapsto\RR, \qquad K(g)=\int_M||\nabla_gJ_g||^2\dif\nu_g,
$$
for which K$\ddot{{\rm a}}$hler metrics are precisely the global
minima. D. Blair and S. Ianus proved that, curiously enough, both
functionals $S$ and $K$ have the same critical points on $\cca$.

\begin{theorem}\label{critherm} \cite{BlrIns}
Let $(M,\omega)$ be a compact symplectic manifold and $\cca$ the set
of all compatible metrics.  Then $g\in\cca$ is a critical point of
$S:\cca\mapsto\RR$ or $K:\cca\mapsto\RR$ if and only if $g$ has
hermitian Ricci tensor.
\end{theorem}

This result and the above considerations do suggest that the
compatible metrics with hermitian Ricci tensor (if any) are really
`good friends' of the symplectic structure.

In \cite{LeWng}, H-V Le and G. Wang approach the problem of the
existence of such metrics by considering an evolution flow inspired
in the Ricci flow introduced by R. Hamilton \cite{Hml1}.  If
$\ricci_g$ is the Ricci tensor of a compatible metric $g$, then
consider the orthogonal decomposition
\begin{equation}\label{ac-c}
\ricci_g=\ricciac_g+\riccic_g,
\end{equation}
where $\ricciac_g=\unm(\ricci_g-\ricci_g(J_g\cdot,J_g\cdot))$ and
$\riccic_g=\unm(\ricci_g+\ricci_g(J_g\cdot,J_g\cdot))$ are the {\it
anti-complexified} and {\it complexified} parts of $\ricci_g$,
respectively. In this way, $g$ has hermitian Ricci tensor if and
only if $\ricciac_g=0$, and since the gradient of the functional $K$
equals $-\ricciac_g$ it is natural to consider the negative gradient
flow equation
\begin{equation}\label{acrf}
\ddt g(t)=\ricciac_{g(t)},
\end{equation}
for a curve $g(t)$ of metrics, which is called in \cite{LeWng} the
{\it anti-complexified Ricci flow}.  Recall that the fixed points of
(\ref{acrf}) are precisely the metrics with hermitian Ricci tensor.
The main result in \cite{LeWng} is the short time existence and
uniqueness of the solution to (\ref{acrf}) when $M$ is compact.

\vs

Let $N$ be a real $2n$-dimensional nilpotent Lie group with Lie
algebra $\ngo$, whose Lie bracket is denoted by $\mu
:\ngo\times\ngo\mapsto\ngo$.  An invariant {\it symplectic}
structure on $N$ is defined by a $2$-form $\omega$ on $\ngo$
satisfying
$$
\omega(X,\cdot)\equiv 0 \quad \mbox{if and only if} \quad X=0 \quad
(\mbox{non-degenerate}),
$$
and for all $X,Y,Z\in\ngo$
\begin{equation}\label{closed}
\omega(\mu(X,Y),Z)+\omega(\mu(Y,Z),X)+\omega(\mu(Z,X),Y)=0 \quad
({\rm closed}, \dif\omega=0).
\end{equation}
Fix a symplectic nilpotent Lie group $(N,\omega)$.  A left invariant
Riemannian metric which is compatible with $(N,\omega)$ is
determined by an inner product $\ip$ on $\ngo$ such that if
$$
\omega(X,Y)=\la X,J_{\ip}Y\ra\quad\forall\; X,Y\in\ngo\quad{\rm
then}\quad J_{\ip}^2=-I.
$$
For the geometric structure $\gamma=\omega$ we have that
$$
G_{\gamma}=\Spe(n,\RR)=\{ g\in\Gl(2n):g^tJg=J\}, \qquad
K_{\gamma}=\U(n),
$$
and the Cartan decomposition of $\ggo_{\gamma}=\spg(n,\RR)=\{
A\in\glg(2n):A^tJ+JA=0\}$ is given by
$$
\spg(n,\RR)=\ug(n)\oplus\pg_{\gamma}, \qquad \pg_{\gamma}=\{
A\in\pg:AJ=-JA\}.
$$
Thus the invariant Ricci tensor $\riccig$ coincides with the
anti-complexified Ricci tensor (see \cite{LeWng}) and for any
$\ip\in\cca$,
\begin{equation}\label{ac-cop}
\Ricg_{\ip}=\Ricac_{\ip}=\unm\left(\Ric_{\ip}+J_{\ip}\Ric_{\ip}J_{\ip}\right).
\end{equation}
This implies that our `goal' condition $\Ricg_{\ip}=0$ (see
(\ref{unconditionint})) is equivalent to have hermitian Ricci
tensor. Also, the evolution flow (\ref{flow}) is not other than the
anti-complexified Ricci flow.

We now review the variational approach developed in Section
\ref{var}.  Fix a non-degenerate $2$-form $\omega$ on $\ngo$, and
let $\Spe(n,\RR)$ denote the subgroup of $\Gl(2n)$ preserving
$\omega$, that is,
$$
\Spe(n,\RR)=\{ \vp\in\Gl(2n):\omega(\vp X,\vp Y)=\omega(X,Y)
\quad\forall\; X,Y\in\ngo\}.
$$
Consider the algebraic subvariety $\nca_s:=\nca_{\gamma}\subset\nca$
given by
$$
\nca_s=\{\mu\in\nca:\dif_{\mu}\omega=0\},
$$
that is, those nilpotent Lie brackets for which $\omega$ is closed
(see (\ref{closed})).  By fixing an inner product $\ip$ on $\ngo$
satisfying that
$$
\omega=\la\cdot,J\cdot\ra\qquad{\rm with}\qquad J^2=-I,
$$
(\ref{ideg}) identifies each $\mu\in\nca_s$ with the
almost-K$\ddot{{\rm a}}$hler manifold $(N_{\mu},\omega,\ip,J)$.  The
action of $\Spe(n,\RR)$ on $\nca_s$ has the following
interpretation:  each $\vp\in\Spe(n,\RR)$ determines a Riemannian
isometry which is also a symplectomorphism
$$
(N_{\vp.\mu},\omega,\ip,J)\mapsto
(N_{\mu},\omega,\la\vp\cdot,\vp\cdot\ra,\vp^{-1}J\vp)
$$
by exponentiating the Lie algebra isomorphism
$\vp^{-1}:(\ngo,\vp.\mu)\mapsto(\ngo,\mu)$.

Let $\ngo$ be a $2n$-dimensional vector space with basis $\{
X_1,...,X_{2n}\}$ over $\RR$, and consider the non-degenerate
$2$-form
$$
\omega=\alpha_1\wedge\alpha_{2n}+...+\alpha_n\wedge\alpha_{n+1},
$$
where $\{\alpha_1,...,\alpha_{2n}\}$ is the dual basis of $\{
X_i\}$.  For the compatible inner product $\la
X_i,X_j\ra=\delta_{ij}$ we have that $\omega=\la\cdot,J\cdot\ra$ for
$$
J=\left[\begin{smallmatrix} &&&&&-1\\ &0&&&\cdot&\\ &&&-1&&\\ &&1&&&\\ &\cdot&&&0&\\
1&&&&&\end{smallmatrix}\right].
$$
In all the following examples the symplectic structure will be
$\omega$, the almost-complex structure $J$ and the compatible metric
$\ip$. We will vary Lie brackets and use constantly identification
(\ref{ideg}).

\begin{example}\label{heis}
{\rm Let $\mu_n$ the $2n$-dimensional Lie algebra whose only
non-zero bracket is
$$
\mu_n(X_1,X_2)=X_3,
$$
that is, $\mu_n$ is isomorphic to $\hg_3\oplus\RR^{2n-3}$, where
$\hg_3$ is the $3$-dimensional Heisenberg Lie algebra.  Recall that
$(N_{\mu_2},\omega)$ is precisely the simply connected cover of the
famous Kodaira-Thurston manifold.  It is easy to prove that
$\overline{\Spe(n,\RR).\mu_n}=\Spe(n,\RR).\mu_n\cup\{ 0\}$, and so
the orbit $\Spe(n,\RR).[\mu_n]$ is closed in $\PP V$.  This implies
that the functional $F$ from Theorem \ref{equiv1gint} must attain
its minimum value on $\Spe(n,\RR).[\mu_n]$ and hence there exists a
metric compatible with $(N_{\mu_n},\omega)$ which is minimal.  In
fact, the inner product $\la X_i,X_j\ra=\delta_{ij}$ satisfies
$$
\begin{array}{l}
\Ricac_{\ip}=-\unc \left[\begin{smallmatrix} 1&&&&&&&&\\ &1&&&&&&&\\ &&-1&&&&&&\\ &&&0&&&&&\\ &&&&\cdot&&&&\\
&&&&&0&&&\\ &&&&&&1&&\\ &&&&&&&-1&\\
&&&&&&&&-1\end{smallmatrix}\right] = -\frac{3}{4}I+\unc
\left[\begin{smallmatrix} 2&&&&&&&&\\ &2&&&&&&&\\ &&4&&&&&&\\ &&&3&&&&&\\
&&&&\cdot&&&&\\ &&&&&3&&&\\ &&&&&&2&&\\ &&&&&&&4&\\ &&&&&&&&4\end{smallmatrix}\right], \qquad 2n\geq 8, \\

\Ricac_{\ip}=-\unc \left[\begin{smallmatrix} 1&&&&&\\ &1&&&&\\ &&-1&&&\\ &&&1&&\\ &&&&-1&\\
&&&&&-1\end{smallmatrix}\right] = -\frac{3}{4}I+\unc \left[\begin{smallmatrix} 2&&&&&\\ &2&&&&\\ &&4&&&\\ &&&2&&\\
&&&&4&\\ &&&&&4\end{smallmatrix}\right], \qquad 2n=6, \\

\Ricac_{\ip}=-\unc \left[\begin{smallmatrix} 1&&&\\ &2&&\\ &&-2&\\
&&&-1\end{smallmatrix}\right] = -\frac{5}{4}I+\unc
\left[\begin{smallmatrix} 4&&&\\ &3&&\\ &&7&\\
&&&6\end{smallmatrix}\right], \qquad 2n=4,
\end{array}
$$
and hence $\Ricac_{\ip}\in\RR I+\Der(\mu_n)$ in all the cases.
Moreover, it follows from the closeness of $\Spe(n,\RR).[\mu_n]$
that $F$ must also attain its maximum value, and therefore
$\Spe(n,\RR).[\mu_n]=\U(n).[\mu_n]$ by uniqueness in Theorem
\ref{equiv1gint}.  This implies that there is only one left
invariant metric compatible with $(N_{\mu_n},\omega)$ up to
isometry, often called the Abbena metric in the case $n=2$. }
\end{example}

\begin{example}\label{otra4}
{\rm Consider the $4$-dimensional Lie algebra given by
$$
\lambda(X_1,X_2)=X_3, \qquad \lambda(X_1,X_3)=X_4.
$$
The compatible metric $\la X_i,X_j\ra=\delta_{ij}$ is minimal for
$(N_{\lambda},\omega)$ since
$$
\Ricac_{\ip}=-\unc \left[\begin{smallmatrix} 3&&&\\ &1&&\\ &&-1&\\
&&&-3\end{smallmatrix}\right] =
-\frac{5}{4}I+\unm\left[\begin{smallmatrix} 1&&&\\ &2&&\\ &&3&\\
&&&4\end{smallmatrix}\right]\in\RR I+\Der(\lambda).
$$
  }
\end{example}

It is well-known that $(N_{\mu_2},\omega)$ and
$(N_{\lambda},\omega)$ are the only symplectic nilpotent Lie groups
in dimension $4$, and then the existence of minimal compatible
metrics in the case $2n=4$ follows.

\begin{example}\label{abc}
{\rm Let $\mu=\mu(a,b,c)$ be the $6$-dimensional $2$-step nilpotent
Lie algebra defined by
$$
\mu(X_1,X_2)=aX_4, \quad \mu(X_1,X_3)=bX_5, \quad \mu(X_2,X_3)=cX_6.
$$
It is easy to check that $\mu\in\nca_s$ if and only if $a-b+c=0$.
We can also get from a simple calculation that
$$
\begin{array}{rl}
\Ricac_{\mu}&=-\unc (a^2+b^2+c^2) \left[\begin{smallmatrix} 1&&&&&\\ &1&&&&\\ &&1&&&\\
&&&-1&&\\ &&&&-1&\\ &&&&&-1\end{smallmatrix}\right] \\
&=-\frac{3}{4}(a^2+b^2+c^2)I+\unm (a^2+b^2+c^2)\left[\begin{smallmatrix} 1&&&&&\\ &1&&&&\\ &&1&&&\\
&&&2&&\\ &&&&2&\\ &&&&&2 \end{smallmatrix}\right]\in\RR I+\Der(\mu),
\end{array}
$$
and so the whole family $\{[\mu(a,b,c)]:a-b+c=0\}\subset\PP\nca_s$
consists of critical points of $F$.  We assume that $a^2+b^2+c^2=2$
in order to avoid homothetical changes, which is equivalent to
$\scalar(\mu)=-1$. The Ricci operator on the center $\zg=\la
X_4,X_5,X_6\ra_{\RR}$ is given by
$$
\Ric_{\mu}=\unm\left[\begin{smallmatrix} a^2&&\\ &b^2&\\
&&c^2\end{smallmatrix}\right],
$$
and thus the curve $\{\mu_{st}=\mu(s,s+t,t):s^2+st+t^2=1,\; 0\leq
t\leq\frac{1}{\sqrt{3}}\}$ is pairwise non-isometric.  It then
follows from Corollary \ref{isoiso} that $(N_{\mu_{st}},\omega)$ is
a curve of pairwise non-isomorphic symplectic nilpotent Lie groups.
In terms of the notation in \cite[Table A.1]{Slm}, we have that
$\mu_{01}\simeq (0,0,0,0,12,13)$ and $\mu_{st}\simeq
(0,0,0,12,13,23)$ for any $0<t\leq\frac{1}{\sqrt{3}}$.  We note that
this curve coincides with the curve of pairwise non-isomorphic
symplectic structures denoted by $\omega_1(t)$ in \cite[Theorem 3.1,
18]{KhkGzMdn}, and then this example shows that any symplectic
structure in such a curve admits a compatible metric which is
minimal  }
\end{example}

\section{Complex structures}\label{complex}

Let $N$ be a real $2n$-dimensional nilpotent Lie group with Lie
algebra $\ngo$, whose Lie bracket is denoted by $\mu
:\ngo\times\ngo\mapsto\ngo$.  An invariant {\it almost-complex}
structure on $N$ is defined by a map $J:\ngo\mapsto\ngo$ satisfying
$J^2=-I$.  If in addition $J$ satisfies the integrability condition
\begin{equation}\label{integral}
\mu(JX,JY)=\mu(X,Y)+J\mu(JX,Y)+J\mu(X,JY), \qquad \forall
X,Y\in\ngo,
\end{equation}
then $J$ is said to be a {\it complex} structure.

Fix an almost-complex nilpotent Lie group $(N,J)$.  A left invariant
Riemannian metric which is {\it compatible} with $(N,J)$, also
called an {\it almost-hermitian metric}, is given by an inner
product $\ip$ on $\ngo$ such that
$$
\la JX,JY\ra=\la X,Y\ra \qquad \forall X,Y\in\ngo.
$$
We have for this particular geometric structure $\gamma=J$ that
$$
G_{\gamma}=\Gl(n,\CC)=\{ \vp\in\Gl(2n):\vp J=J\vp\}, \qquad
K_{\gamma}=\U(n),
$$
and the Cartan decomposition of $\ggo_{\gamma}=\glg(n,\CC)=\{
A\in\glg(2n):AJ=JA\}$ is given by
$$
\glg(n,\CC)=\ug(n)\oplus\pg_{\gamma}, \qquad \pg_{\gamma}=\{
A\in\pg:AJ=JA\}.
$$
The invariant Ricci operator is then given by the complexified Ricci
operator
$$
\Ricg_{\ip}=\Ricc_{\ip}=\unm\left(\Ric_{\ip}-J\Ric_{\ip}J\right)
$$
(see (\ref{ac-c})).  In this way, condition $\Ricg_{\ip}=0$ is
equivalent to the Ricci operator anti-commute with $J$. We do not
know if this property has any special significance in complex
geometry, but for instance it holds for a K$\ddot{{\rm a}}$hler
metric if and only if the metric is Ricci flat.  Anyway, as in the
symplectic case, the condition $\Ricc_{\ip}=0$ is also forbidden for
non-abelian $N$ since $\tr{\Ricc_{\ip}}=\scalar(\ip)<0$.

We now fix a map $J:\ngo\mapsto\ngo$ satisfying $J^2=-I$ and
consider the algebraic subvariety $\nca_c:=\nca_{\gamma}\subset\nca$
given by
$$
\nca_c=\{\mu\in\nca: \mbox{(\ref{integral}) holds}\},
$$
that is, those nilpotent Lie brackets for which $J$ is integrable
and so define a complex structure on $N_{\mu}$, the simply connected
nilpotent Lie group with Lie algebra $(\ngo,\mu)$.

Fix also an inner product $\ip$ on $\ngo$ compatible with $J$, then
(\ref{ideg}) identifies each $\mu\in\nca_c$ (or $\nca$) with the
hermitian (or almost-hermitian) manifold $(N_{\mu},J,\ip)$.  If we
use the same triple $(\omega,J,\ip)$ to define and identify $\nca_s$
(see Section \ref{symp}) and $\nca_c$, then the intersection of
these varieties is $\nca_s\cap\nca_c=\{ 0\}$ since no non-abelian
nilpotent Lie group can admit a left invariant K$\ddot{{\rm a}}$hler
metric.

We now give some examples.

\begin{example}\label{heisc}
{\rm Let $\mu_n$ be the $2n$-dimensional Lie algebra considered in
Example \ref{heis}.  It is easy to check that $\ip$ is also minimal
as a compatible metric for the almost-complex nilpotent Lie group
$(N_{\mu_n},J)$.  For the $4$-dimensional Lie algebra in Example
\ref{otra4}, we have that $\Ricc_{\ip}=-\unc I$ and hence this
metric is minimal for the almost-complex nilpotent Lie group
$(N_{\lambda},J)$ as well.  }
\end{example}

For $\ngo_1=\RR^4$ and $\ngo_2=\RR^2$, consider the vector space
$W=\Lambda^2\ngo_1^*\otimes\ngo_2$ of all skew-symmetric bilinear
maps $\mu:\ngo_1\times\ngo_1\mapsto\ngo_2$.  Any $6$-dimensional
$2$-step nilpotent Lie algebra with $\dim{\mu(\ngo,\ngo)}\leq 2$ can
be modelled in this way.  Fix basis $\{ X_1,...,X_4\}$ and $\{
Z_1,Z_2\}$ of $\ngo_1$ and $\ngo_2$, respectively.  Each element in
$W$ will be described as $\mu=\mu(a_1,a_2,...,f_1,f_2)$, where
{\small
$$
\begin{array}{lll}
\mu(X_1,X_2)=a_1Z_1+a_2Z_2, & \mu(X_1,X_4)=c_1Z_1+c_2Z_2, & \mu(X_2,X_4)=e_1Z_1+e_2Z_2,\\
\mu(X_1,X_3)=b_1Z_1+b_2Z_2, & \mu(X_2,X_3)=d_1Z_1+d_2Z_2, &
\mu(X_3,X_4)=f_1Z_1+f_2Z_2.
\end{array}
$$}
The complex structure and the compatible metric will be always
defined by
$$
\begin{array}{lll}
J=\left[\begin{smallmatrix} 0&-1&&&&\\ 1&0&&&&\\ &&0&-1&&\\
&&1&0&&\\ &&&&0&-1\\ &&&&1&0
\end{smallmatrix}\right], && \la X_i,X_j\ra=\la Z_i,Z_j\ra=\delta_{ij}.
\end{array}
$$
If $A=(a_1,a_2)$, ... , $F=(f_1,f_2)$ and $JA=(-a_2,a_1)$, ... ,
$JF=(-f_2,f_1)$, then it is easy to check that $J$ is integrable on
$N_{\mu}$ (or $(N_{\mu},J)$ is a complex nilpotent Lie group) if and
only if
\begin{equation}\label{integrable}
E=B+JD+JC,
\end{equation}
$J$ is {\it bi-invariant} (i.e. $\mu(JX,Y)=J\mu(X,Y)$) if and only
if
\begin{equation}\label{biinv}
A=F=0, \qquad C=D=JB, \qquad E=-B,
\end{equation}
and $J$ is {\it abelian} (i.e. $\mu(JX,JY)=\mu(X,Y)$) if and only if
\begin{equation}\label{abelian}
E=B, \qquad D=-C.
\end{equation}
We note that the above conditions determine
$\Gl(2,\CC)\times\Gl(1,\CC)$-invariant linear subspaces of $W$ of
dimensions $10$, $2$ and $8$, respectively.  For each $\mu\in W$, it
is easy to see that the Ricci operator of the almost-hermitian
manifold $(N_{\mu},J,\ip)$ (see identification (\ref{ideg}))
restricted to $\ngo_1$, $\Ric_{\mu}|_{\ngo_1}$, is given by
\begin{equation}\label{ricciw}
-\unm\left[\begin{smallmatrix}
||A||^2+||B||^2+||C||^2 & \la B,D\ra+\la C,E\ra & -\la A,D\ra+\la C,F\ra & -\la A,E\ra-\la B,F\ra\\
\la B,D\ra+\la C,E\ra & ||A||^2+||D||^2+||E||^2 & \la A,B\ra+\la E,F\ra & \la A,C\ra-\la D,F\ra\\
-\la A,D\ra+\la C,F\ra & \la A,B\ra+\la E,F\ra & ||B||^2+||D||^2+||F||^2 & \la B,C\ra+\la D,E\ra\\
-\la A,E\ra-\la B,F\ra & \la A,C\ra-\la D,F\ra & \la B,C\ra+\la
D,E\ra & ||C||^2+||E||^2+||F||^2
\end{smallmatrix}\right]
\end{equation}
and
$$
\Ric_{\mu}|_{\ngo_2}=\unm\left[\begin{smallmatrix} ||v_1||^2 & \la
v_1,v_2\ra\\ \la v_1,v_2\ra & ||v_2||^2
\end{smallmatrix}\right], \quad v_i=(a_i,b_i,c_i,d_i,e_i,f_i),\quad i=1,2.
$$
Recall that if the complexified Ricci operator satisfied
$\Ricc_{\mu}|_{\ngo_1}=pI$, $p\in\RR$, then $\mu$ is minimal.
Indeed, since always $\Ricc_{\mu}|_{\ngo_2}=qI$ for some $q\in\RR$,
we would have that
\begin{equation}\label{multide}
\Ricc_{\mu}=\left[\begin{smallmatrix}pI&\\&qI\end{smallmatrix}\right]=(2p-q)I+
\left[\begin{smallmatrix}(q-p)I&\\
&2(q-p)I\end{smallmatrix}\right]\in\RR I+\Der(\mu).
\end{equation}
In particular, any bi-invariant complex nilpotent Lie group
$(N_{\mu},J)$ (see (\ref{biinv})) admits a compatible metric which
is minimal.

We will now focus on the abelian complex case (see (\ref{abelian})).
It is not hard to see that these conditions imply that
$\Ricc_{\mu}|_{\ngo_1}=\Ric_{\mu}$, and so in this case, to get
$\Ricc_{\mu}|_{\ngo_1}\in\RR I$ is necessary and sufficient that
$$
\la A+F,B\ra=0, \qquad \la A+F,C\ra=0, \qquad ||A||^2=||F||^2.
$$
In order to avoid homothetical changes we will always ask for
$||v_1||^2+||v_2||^2=2$, which is equivalent to $\scalar(\mu)=-1$.

\begin{example} {\rm If we put $A=(s,t)$, $F=(-s,t)$, $B=C=D=E=0$, $s^2+t^2=1$, then the corresponding curve
$\mu_{st}$ of minimal compatible metrics satisfies
$$
\Ric_{\mu_{st}}|_{\ngo_2}=\left[\begin{smallmatrix} s^2&0\\
0&t^2\end{smallmatrix}\right],
$$
proving that $\{\mu_{st}:s^2+t^2=1,\; 0\leq s\leq
\frac{1}{\sqrt{2}}\}$ is a curve of pairwise non-isometric metrics.
It then follows from Corollary \ref{isoiso} that $(N_{\mu_{st}},J)$
is a curve of pairwise non-isomorphic abelian complex nilpotent Lie
groups.  Recall that $\mu_{st}\simeq\hg_3\oplus\hg_3$ for all $0<s$
and $\mu_{01}\simeq\hg_3\oplus\RR^3$.  }
\end{example}

\begin{example} {\rm For $A=(s,t)$, $F=(-s,t)$, $B=(\unm,0)=C=-D=E$, $s^2+t^2=\unm$, the curve
$\mu_{st}$ of minimal compatible metrics satisfies
$$
\Ric_{\mu_{st}}|_{\ngo_2}=\left[\begin{smallmatrix} s^2+\unm&0\\
0&t^2\end{smallmatrix}\right],
$$
which implies that the family  $\{\mu_{st}:s^2+t^2=\unm\}$ is
pairwise non-isometric. It is easy to see that for $t\ne 0$,
$\mu_{st}$ is isomorphic to the complex Heisenberg Lie algebra, and
hence $(N_{\mu_{st}},J)$ defines a curve of pairwise non-isomorphic
abelian complex structures on the Iwasawa manifold. Since
$j_{\mu_{st}}(Z_2)^2\notin\RR I$, we have that the hermitian
manifolds $(N_{\mu_{st}},J,\ip)$ are not modified H-type (see
\cite{modified}).  }
\end{example}

\begin{example} {\rm Consider the abelian complex structures defined by $A=-F$, $E=B$ and $D=-C$.  In this case,
the Hermitian manifolds $(N_{\mu},J,\ip)$ are modified H-type and
$\mu$ is always isomorphic to the complex Heisenberg Lie algebra
when $v_1,v_2\ne 0$.  In fact, by assuming for simplicity that $\la
v_1,v_2\ra=0$, then
$$
j_{\mu}(Z)^2=-\unm(\la Z,Z_1\ra^2 ||v_1||^2+\la Z,Z_2\ra^2
||v_2||^2)I, \qquad \forall Z\in\ngo_2.
$$
For $A=(s,0)=-F$, $B=(0,t)=E$, $s^2+t^2=1$, $D=C=0$, the
corresponding curve $\mu_{st}$ of minimal compatible metrics
satisfies
$$
\Ric_{\mu_{st}}|_{\ngo_2}=\left[\begin{smallmatrix} s^2&0\\
0&t^2\end{smallmatrix}\right],
$$
and so the family  $\{\mu_{st}:s^2+t^2=1,\; 0\leq s\leq
\frac{1}{\sqrt{2}} \}$ is pairwise non-isometric and the abelian
complex structures $(N_{\mu_{st}},J)$ are pairwise non-isomorphic.
Each modified H-type metric is compatible with two spheres of
abelian complex structures of this type which can be described by
$$
\{ \pm v_1\times v_2:v_i\in\RR^3,\; ||v_1||^2=2s^2,\;
||v_2||^2=2t^2,\; \la v_1,v_2\ra=0\},
$$
where $v_1\times v_2$ denotes the vectorial product, but one can see
that these structures are all isomorphic to $(N_{\mu_{st}},J)$
(compare with \cite{KtsSlm}).  We finally recall that
$\mu_{01}\simeq\hg_5\oplus\RR$, and so $\ip$ is a minimal compatible
metric for the abelian complex nilpotent Lie group
$(N_{\mu_{01}},J)$. }
\end{example}

Although it has not been mentioned, most of the curves given in this
section have been obtained via the variational method provided by
Theorem \ref{equiv1gint}, by using an approach very similar to that
in \cite[Example 5.2]{minimal}.

\section{Hypercomplex structures}\label{hyper}

Let $N$ be a real $4n$-dimensional nilpotent Lie group with Lie
algebra $\ngo$, whose Lie bracket is denoted by $\mu
:\ngo\times\ngo\mapsto\ngo$.  An invariant {\it hypercomplex}
structure on $N$ is defined by a triple $\{ J_1,J_2,J_3\}$ of
complex structures on $\ngo$ (see Section \ref{complex}) satisfying
the quaternion identities
\begin{equation}\label{quat}
J_i^2=-I,\quad i=1,2,3, \qquad J_1J_2=J_3=-J_2J_1.
\end{equation}
An inner product $\ip$ on $\ngo$ is said to be {\it compatible} with
$\{ J_1,J_2,J_3\}$, also called an {\it hyper-hermitian metric}, if
\begin{equation}\label{ortconhyp}
\la J_iX,J_iY\ra=\la X,Y\ra \qquad \forall X,Y\in\ngo, \; i=1,2,3.
\end{equation}
Two hypercomplex nilpotent Lie groups $(N,\{ J_1,J_2,J_3\})$ and
$(N',\{ J_1',J_2',J_3'\})$ are said to be {\it isomorphic} if there
exists an automorphism $\vp:\ngo'\mapsto\ngo$ such that
$$
\vp J_i'\vp^{-1}=J_i, \quad i=1,2,3.
$$
For $\gamma=\{ J_1,J_2,J_3\}$ we therefore have that
$$
G_{\gamma}=\Gl(n,\HH)=\{ \vp\in\Gl(4n):\vp J_i=J_i\vp,\; i=1,2,3\},
\qquad K_{\gamma}=\Spe(n),
$$
and the Cartan decomposition of
$$
\ggo_{\gamma}=\glg(n,\HH)=\{ A\in\glg(4n):AJ_i=J_iA,\; i=1,2,3\}
$$
is given by
$$
\glg(n,\HH)=\spg(n)\oplus\pg_{\gamma}, \qquad \pg_{\gamma}=\{
A\in\pg:AJ_i=J_iA,\; i=1,2,3\}.
$$
The invariant Ricci operator for a compatible metric $\ip\in\cca$ is
then given by
$$
\Ricg_{\ip}=\unc(\Ric_{\ip}-J_1\Ric_{\ip}J_1-J_2\Ric_{\ip}J_2-J_3\Ric_{\ip}J_3),
$$
and hence condition $\Ricg_{\ip}=0$ can never holds since
$\tr{\Ricg_{\ip}}=\scalar(\ip)<0$ for a non-abelian nilpotent Lie
group.

We now fix $\{ J_1,J_2,J_3,\ip\}$ satisfying (\ref{quat}) and
(\ref{ortconhyp}) and consider the algebraic subvariety
$\nca_h:=\nca_{\gamma}\subset\nca$ given by
$$
\nca_h=\{\mu\in\nca: J_i \;\mbox{is integrable for}\; i=1,2,3\},
$$
that is, those nilpotent Lie brackets for which $\{ J_1,J_2,J_3\}$
is an hypercomplex structure on the corresponding nilpotent Lie
group $N_{\mu}$.  Thus (\ref{ideg}) identifies each $\mu\in\nca_h$
with the hyper-hermitian manifold $(N_{\mu},\{ J_1,J_2,J_3\},\ip)$.

It is proved in \cite{minimal} that any hypercomplex $8$-dimensional
nilpotent Lie group admits a minimal compatible metric, and that the
moduli space of all hypercomplex $8$-dimensional nilpotent Lie
groups up to isomorphism is $9$-dimensional, and the moduli space of
the abelian ones has dimension $5$.

We will give now explicit continuous families of hypercomplex
structures on some particular nilpotent Lie groups.

For $\ngo_1=\RR^4$ and $\ngo_2=\RR^4$ consider the vector space
$W=\Lambda^2\ngo_1^*\otimes\ngo_2$ of all skew-symmetric bilinear
maps $\mu:\ngo_1\times\ngo_1\mapsto\ngo_2$.  Any $8$-dimensional
$2$-step nilpotent Lie algebra with $\dim{\mu(\ngo,\ngo)}\leq 4$ can
be modelled in this way.  Fix basis $\{ X_1,X_2,X_3,X_4\}$ and $\{
Z_1,Z_2,Z_3,Z_4\}$ of $\ngo_1$ and $\ngo_2$, respectively.  Each
element in $W$ will be denoted as
$\mu=\mu(a_1,...,a_4,...,f_1,...,f_4)$, where {\small
$$
\begin{array}{ll}
\mu(X_1,X_2)=a_1Z_1+a_2Z_2+a_3Z_3+a_4Z_4, & \mu(X_2,X_3)=d_1Z_1+d_2Z_2+d_3Z_3+d_4Z_4,\\
\mu(X_1,X_3)=b_1Z_1+b_2Z_2+b_3Z_3+b_4Z_4, & \mu(X_2,X_4)=e_1Z_1+e_2Z_2+e_3Z_3+e_4Z_4, \\
\mu(X_1,X_4)=c_1Z_1+c_2Z_2+c_3Z_3+c_4Z_4, &
\mu(X_3,X_4)=f_1Z_1+f_2Z_2+f_3Z_3+f_4Z_4.
\end{array}
$$ }
The compatible metric will be $\la X_i,X_j\ra=\la
Z_i,Z_j\ra=\delta_{ij}$ and the hypercomplex structure will always
act on $\ngo_i$ by
$$
J_1=\left[\begin{smallmatrix} 0&-1&&\\ 1&0&&\\ &&0&-1\\
&&1&0\end{smallmatrix}\right], \quad J_2=\left[\begin{smallmatrix}
&&-1&0\\ &&0&1\\ 1&0&&\\ 0&-1&&\end{smallmatrix}\right],\quad
J_3=\left[\begin{smallmatrix} &&&-1\\ &&-1&\\ &1&&\\
1&&&\end{smallmatrix}\right].
$$
If $A=(a_1,...,a_4)$, ..., $F=(f_1,...,f_4)$, then it is easy to
prove that $J_i$ is integrable for all $i=1,2,3$ on $N_{\mu}$ (or
$(N_{\mu},\{ J_1,J_2,J_3\})$ is a hypercomplex nilpotent Lie group)
if and only if
$$
E=B+J_1D+J_1C, \qquad D=-C-J_2A-J_2F, \qquad F=-A-J_3B+J_3E.
$$
If we define $T:=D+C$, then the above conditions are equivalent to
\begin{equation}\label{integrable3}
D=-C+T, \qquad E=B+J_1T, \qquad F=-A+J_2T.
\end{equation}
In order to use a notation as similar as possible to \cite{DttFin1,
DttFin2}, we should put $T=(t_3,t_2,-t_1,t_4)$.  It is easy to check
that $(N_{\mu},\{ J_1,J_2,J_3\})$ is {\it abelian} (i.e.
$\mu(J_i\cdot,J_i\cdot)=\mu$, $i=1,2,3$) if and only if $T=0$.  We
note that $\dim{W}=24$, and so condition (\ref{integrable3})
determine a $\Gl(1,\HH)\times\Gl(1,\HH)$-invariant linear subspace
$W_h$ of $W$ of dimension $16$, and a $12$-dimensional subspace
$W_{ah}$ if we ask in addition abelian.

For each $\mu\in W$, the Ricci operator of $(N_{\mu},\{
J_1,J_2,J_3\},\ip)$ (see identification (\ref{ideg})) restricted to
$\ngo_1$, $\Ric_{\mu}|_{\ngo_1}$, is given by formula
(\ref{ricciw}), and
$$
\Ric_{\mu}|_{\ngo_2}=\unm[\la v_i,v_j\ra], \quad 1\leq i,j\leq 4,
\quad v_i:=(a_i,b_i,c_i,d_i,e_i,f_i).
$$
Since the only symmetric transformations of $\ngo_i=\RR^4$ commuting
with all the $J_i's$ are the multiplies of the identity, we obtain
that the invariant Ricci operator satisfies
$\Ricg_{\mu}|_{\ngo_i}\in\RR I$ for any $\mu\in W$.  By arguing as
in (\ref{multide}), one obtains that any $\mu\in W$ is minimal.

Let $\ggo_1$, $\ggo_2$ and $\ggo_3$ denote the $8$-dimensional Lie
algebras obtained as the direct sum of an abelian factor and the
following H-type Lie algebras: the $5$-dimensional Heisenberg Lie
algebra, the $6$-dimensional complex Heisenberg Lie algebra and the
$7$-dimensional quaternionic Heisenberg Lie algebra.  In order to
avoid homothetical changes we will always ask for
$||v_1||^2+...+||v_4||^2=2$, which is equivalent to
$\scalar(\mu)=-1$.

\begin{example} {\rm If we put $T=0$, $A=(0,r,0,0)$, $B=(0,0,s,0)$, $C=(0,0,0,t)$,
we have for
each $\mu_{rst}$ that
$$
\Ric_{\mu_{rst}}|_{\ngo_2}=\unm\left[\begin{smallmatrix} 0&&&\\ &r^2&&\\ &&s^2&\\
&&&t^2\end{smallmatrix}\right],
$$
and thus the family
$$
\{ (N_{\mu_{rst}},\{ J_1,J_2,J_3\},\ip):0\leq r\leq s\leq t,\;
r^2+s^2+t^2=2\}
$$
of minimal compatible metrics is pairwise non-isometric.  This gives
rise then a $2$-parameter family of pairwise non-isomorphic abelian
hypercomplex nilpotent Lie groups (see Corollary \ref{isoiso}).  If
$0<r$ then $\mu_{rst}\simeq\ggo_3$, for $r=0<s$ we get a curve on
$\ggo_2$ and for $r=s=0$, $t=1$, a single structure on $\ggo_1$.   }
\end{example}

\begin{example}\label{5g3} {\rm We now set $T=0$ and choose $A,B,C$ such that {\small
$$
v_1=(\displaystyle{\frac{1}{\sqrt{2}}},0,0,0,0,-\frac{1}{\sqrt{2}}),
\; v_2=(0,\sqrt{\frac{3}{8}},0,0,\sqrt{\frac{3}{8}},0), \;
||v_3||^2+||v_4||^2=\frac{1}{4}, \; ||v_3||^2>||v_4||^2.
$$ }
Assume that two of such elements $\lambda=\mu(v_3,v_4)$ and
$\lambda'=\mu(v_3',v_4')$ are in the same
$\Spe(1)\times\Spe(1)$-orbit, say $\lambda'=\vp.\lambda$ for some
$\vp=(\vp_1,\vp_2)\in\Spe(1)\times\Spe(1)$. Recall that
$$
j_{\lambda'}(Z)=\vp_1j_{\lambda}(\vp_2^{-1}Z)\vp_1^{-1}, \qquad
\forall Z\in\ngo_2,
$$
(see \cite[Appendix]{minimal}) and $\la
v_i,v_j\ra=-\unm\tr{j_{\mu}(Z_i)j_{\mu}(Z_j)}$, $1\leq i,j\leq 4$.
It then follows from
$$
||v_1||>||v_2||>||v_3||>||v_4||, \qquad
||v_1'||>||v_2'||>||v_3'||>||v_4'||
$$
that $j_{\vp.\lambda}(Z_i)=\pm j_{\lambda}(Z_i)$ for all
$i=1,...,4$, and hence $v_3'=\pm v_3$ and $v_4'=\pm v_4$.  Thus we
have a family of pairwise non-isomorphic abelian hypercomplex
structures on $\ggo_3$ depending on $5$ parameters (see Corollary
\ref{isoiso}).  Analogously, we get a $5$-dimensional family on
$\ggo_2$ by putting $v_1=v_2=0$. }
\end{example}

\begin{example}
{\rm Let $\mu_t$ be the curve defined for $0\leq
t\leq\frac{1}{\sqrt{3}}$ by {\small
$$
\begin{array}{lll}
\mu_{t}(X_1,X_2)=\sqrt{1-3t^3}Z_1+tZ_2, && \mu_{t}(X_2,X_3)=tZ_4, \\
\mu_{t}(X_1,X_3)=tZ_3, && \mu_{t}(X_2,X_4)=-tZ_3, \\
\mu_{t}(X_1,X_4)=tZ_4, && \mu_{t}(X_3,X_4)=-\sqrt{1-3t^3}Z_1+tZ_2.
\end{array}
$$ }
It is easy to check that $(N_{\mu_t},\{ J_1,J_2,J_3\})$ is always
non-abelian hypercomplex (recall that $t_1=t_2=t_3=0$, $t_4=2t$) and
the curve is pairwise non-isomorphic since it follows from
$$
\Ric_{\mu_{t}}|_{\ngo_2}=\left[\begin{smallmatrix} 1-3t^2&&&\\ &t^2&&\\ &&t^2&\\
&&&t^2\end{smallmatrix}\right]
$$
that the curve $(N_{\mu_t},\ip)$ is pairwise non-isometric (see
Corollary \ref{isoiso}).  The starting and ending points are
$\mu_0\simeq\ggo_1$ and $\mu_{\frac{1}{\sqrt{3}}}\simeq\ggo_3$,
respectively, and $\mu_t\simeq\ug(2)\oplus\CC^2$ for any
$0<t<\frac{1}{\sqrt{3}}$ (see \cite{manus} for further information
on these $2$-step nilpotent Lie algebras constructed via
representations of compact Lie groups).  }
\end{example}

\section{Einstein solvmanifolds}\label{einstein}

Our goal in this section is to show that the `moment map' approach
proposed in \cite{minimal} can also be applied to the study of
Einstein solvmanifolds.  After a brief overview of such spaces, we
will follow the same path used to study compatible metrics in the
previous sections, but by considering the Ricci operator
$\Ric_{\ip}$ itself. In other words, none geometric structure is
considered (or $\gamma=0$).

A {\it solvmanifold} is a solvable Lie group $S$ endowed with a left
invariant Riemannian metric, and $S$ is called {\it standard} if
$\ag:=\ngo^{\perp}$ is abelian, where $\ngo=[\sg,\sg]$ and $\sg$ is
the Lie algebra of $S$.  All known examples of non-compact
homogeneous Einstein manifolds are isometric to standard Einstein
solvmanifolds.  These spaces have been extensively studied by J.
Heber in \cite{Hbr}, obtaining remarkable structure and uniqueness
results.

Let $N$ be a nilpotent Lie group with Lie algebra $\ngo$ of
dimension $n$, whose Lie bracket is denoted by $\mu
:\ngo\times\ngo\mapsto\ngo$.  We have in this case $\gamma=0$, thus
any inner product is `compatible', $G_{\gamma}=\Gl(n)$,
$K_{\gamma}=\Or(n)$, $\pg_{\gamma}=\pg$, $\nca_{\gamma}=\nca$,
$\Ricg=\Ric$ and then condition $\Ricg=0$ is clearly forbidden for
non-abelian $N$.  Moreover, the evolution equation is precisely the
Ricci flow and the corresponding invariant Ricci solitons coincide
with minimal metrics.

Given a metric nilpotent Lie algebra $(\ngo,\ip)$, a metric solvable
Lie algebra $(\sg=\ag\oplus\ngo,\ip')$ is called a {\it metric
solvable extension} of $(\ngo,\ip)$ if the restrictions of the Lie
bracket of $\sg$ and the inner product $\ip'$ to $\ngo$ coincide
with the Lie bracket of $\ngo$ and $\ip$, respectively.  It turns
out that for each $(\ngo,\ip)$ there exists a unique rank-one (i.e.
$\dim{\ag}=1$) metric solvable extension of $(\ngo,\ip)$ which stand
a chance of being an Einstein space (see \cite{critical}).  This
fact turns the study of rank-one Einstein solvmanifolds into a
problem on nilpotent Lie algebras.  More specifically, a nilpotent
Lie algebra $\ngo$ is the nilradical of a standard Einstein
solvmanifold if and only if $\ngo$ admits an inner product $\ip$
such that $\Ric_{\ip}=cI+D$ for some $c\in\RR$ and $D\in\Der(\ngo)$,
that is, $\ip$ is minimal.

Let $\nca$ be the variety of all nilpotent Lie algebras of dimension
$n$.  Fix an inner product $\ip$ on $\ngo$. Each $\mu\in\nca$ is
then identified via (\ref{ideg}) with the Riemannian manifold
$(N_{\mu},\ip)$, but we also have in this case another
identification with a solvmanifold: for each $\mu\in\nca$, there
exists a unique rank-one metric solvable extension
$S_{\mu}=(S_{\mu},\ip)$ of $(N_{\mu},\ip)$ standing a chance of
being Einstein, and every $(n+1)$-dimensional rank-one Einstein
solvmanifold can be modelled as $S_{\mu}$ for a suitable
$\mu\in\nca$.  We recall that the study of standard solvmanifolds
reduces essentially to the rank-one case (see \cite[4.18]{Hbr}).

The functional $F:\PP \nca\mapsto\RR$ given by
$F([\mu])=\tr{\Ric_{\mu}^2}/||\mu||^4$ measures how far is the
metric $\mu$ from being Einstein (see \cite{soliton}).

From Theorem \ref{equiv1gint} we then obtain the uniqueness up to
isometry of Einstein metrics on standard solvable Lie groups proved
in \cite{Hbr}, as well as the variational result given in
\cite{critical} characterizing Einstein solvmanifolds as critical
points of a natural curvature functional.

Theorem \ref{equiv2gint} gives the relationship between Ricci
soliton metrics on nilpotent Lie groups and Einstein solvmanifolds
proved in \cite{soliton}.  Part (i) is a new characterization of
these privileged metrics, claiming that they are precisely minimal
metrics.

As far as we know, the following is a complete list of the known
examples of nilpotent Lie groups admitting a minimal metric, or
equivalently, of the nilradicals of standard Einstein solvmanifolds:
\begin{itemize}

\item Iwasawa $N$-groups: $G/K$ irreducible symmetric space of
noncompact type and $G=KAN$ the Iwasawa decomposition.

\item\cite{PttS} Nilradicals of normal $j$-algebras (i.e. of
noncompact homogeneous K$\ddot{{\rm a}}$hler Einstein spaces).

\item\cite{Dlf} Certain $2$-step nilpotent Lie algebras for which there
is a basis with very uniform properties (see also \cite[1.9]{Wlt}).

\item\cite{Bgg} $H$-type Lie groups (see also \cite{Lnz}).

\item\cite{Crt} Nilradicals of homogeneous quaternionic K$\ddot{{\rm
a}}$hler spaces.

\item\cite{EbrHbr, manus} Naturally reductive nilpotent Lie groups.

\item\cite{Hbr} Families of deformations of Iwasawa $N$-groups in the rank-one
case.

\item\cite{Fan1,Fan2} Certain $2$-step nilpotent Lie algebras
constructed via Clifford modules.

\item\cite{GrdKrr} A $2$-parameter family of $9$-dimensional
$2$-step nilpotent Lie algebras (with $3$-dimensional center) and
certain modifications of Iwasawa $N$-groups (rank $\geq 2$).

\item\cite{finding} Nilpotent Lie algebras with a codimension one abelian ideal.

\item\cite{finding} A curve of $6$-step nilpotent Lie algebras of dimension $7$, which
is the lowest possible dimension for a continuous family.

\item\cite{Ymd, Mor} Certain $2$-step nilpotent Lie algebras defined from
subsets of fundamental roots of complex simple Lie algebras.

\item\cite{Wll, finding} Nilpotent Lie algebras of dimension $\leq 6$.

\item\cite{inter} A curve of $10$-dimensional $2$-step nilpotent Lie
algebras with $5$-dimensional center.

\item\cite{Krr} A $2$-parameter family of deformations of the
$12$-dimensional quaternionic hyperbolic space.

\item\cite{Tmr} Nilradicals of parabolic subalgebras of semisimple
Lie algebras which are $2$-step or $3$-step.

\end{itemize}

\end{document}